\documentclass[11pt]{article}
\usepackage{mathrsfs}
\usepackage{amssymb}
\usepackage{amsmath}
\usepackage{amsbsy}
\usepackage{epsfig}
\usepackage{enumerate}
\usepackage{bm}
\usepackage{color}

\usepackage[colorlinks,linkcolor=blue]{hyperref}

\newtheorem{lemma}{Lemma}[section]
\newtheorem{theorem}{Theorem}[section]
\newtheorem{proposition}{Proposition}[section]

\newtheorem{corollary}{Corollary}[section]
\newtheorem{definition}{Definition}[section]
\newtheorem{remark}{Remark}[section]

\def\pr{\textsf{P}} 
\def\ep{\textsf{E}} 
\def\Sbep{\widehat{\mathbb E}} 
\def\cSbep{\widehat{\mathcal E}} 
\def\Capc{\mathbb V} 
\def\cCapc{\mathcal V} 
\def\upCapc{\widehat{\mathbb V}} 
\def\lowCapc{\widehat{\mathcal V}} 
\def\outCapc{\widehat{\mathbb V}^{\ast}}
\def\outcCapc{\widehat{\mathcal V}^{\ast}} 

\def\tildeCapc{\mathbb C^{\ast}}

\def\vSbep{\breve{\mathbb E}}

\renewcommand{\baselinestretch}{1.5}

\begin{document}

\begin{center}{\LARGE\bf A note on the cluster set of the law  of the iterated logarithm   under sub-linear expectations}
\end{center}

\begin{center} {\sc
Li-Xin Zhang\footnote{This work was supported by grants from the NSF of China (Nos. 11731012,12031005),   Ten Thousands Talents Plan of Zhejiang Province (Grant No. 2018R52042), NSF of Zhejiang Province (No. LZ21A010002) and the Fundamental Research Funds for the Central Universities
}
}\\
{\sl \small School  of Mathematical Sciences, Zhejiang University, Hangzhou 310027} \\
(Email:stazlx@zju.edu.cn)\\
\end{center}

 \renewcommand{\abstractname}{~}
\begin{abstract}
{\bf Abstract:} In this note, we establish   a compact law of the iterated logarithm under the upper capacity  for independent and identically distributed random variables   in a sub-linear expectation space. For showing the result, a self-normalized law of the iterated logarithm is established.




\vspace{-3mm}
\end{abstract}

\renewcommand{\baselinestretch}{1.2}


\section{ Introduction and notations.}\label{sect1}
\setcounter{equation}{0}

We use the framework and notations of Peng (2008,   2019). Let  $(\Omega,\mathcal F)$
 be a given measurable space  and let $\mathscr{H}$ be a linear space of real measurable functions
defined on $(\Omega,\mathcal F)$ such that if $X_1,\ldots, X_n \in \mathscr{H}$  then $\varphi(X_1,\ldots,X_n)\in \mathscr{H}$ for each
$\varphi\in C_{l,Lip}(\mathbb R^n)$,  where $C_{l,Lip}(\mathbb R^n)$ denotes the linear space of (local Lipschitz)
functions $\varphi$ satisfying
\begin{eqnarray*} & |\varphi(\bm x) - \varphi(\bm y)| \le  C(1 + |\bm x|^m + |\bm y|^m)|\bm x- \bm y|, \;\; \forall \bm x, \bm y \in \mathbb R^n,&\\
& \text {for some }  C > 0, m \in \mathbb  N \text{ depending on } \varphi. &
\end{eqnarray*}
$\mathscr{H}$ is considered as a space of ``random variables''.   We also denote $C_{b,Lip}(\mathbb R^n)$ the space of bounded  Lipschitz
functions.

\begin{definition}\label{def1.1} A  sub-linear expectation $\Sbep$ on $\mathscr{H}$  is a function $\Sbep: \mathscr{H}\to \overline{\mathbb R}$ satisfying the following properties: for all $X, Y \in \mathscr H$, we have
\begin{description}
  \item[\rm (a)]  Monotonicity: If $X \ge  Y$ then $\Sbep [X]\ge \Sbep [Y]$;
\item[\rm (b)] Constant preserving: $\Sbep [c] = c$;
\item[\rm (c)] Sub-additivity: $\Sbep[X+Y]\le \Sbep [X] +\Sbep [Y ]$ whenever $\Sbep [X] +\Sbep [Y ]$ is not of the form $+\infty-\infty$ or $-\infty+\infty$;
\item[\rm (d)] Positive homogeneity: $\Sbep [\lambda X] = \lambda \Sbep  [X]$, $\lambda\ge 0$.
 \end{description}
 Here $\overline{\mathbb R}=[-\infty, \infty]$, $0\cdot \infty$ is defined to be $0$. The triple $(\Omega, \mathscr{H}, \Sbep)$ is called a sub-linear expectation space. Given a sub-linear expectation $\Sbep $, let us denote the conjugate expectation $\cSbep$of $\Sbep$ by
$$ \cSbep[X]:=-\Sbep[-X], \;\; \forall X\in \mathscr{H}. $$
\end{definition}
 By Theorem 1.2.1 of Peng (2019), there exists a family of finite additive linear expectations $E_{\theta}: \mathscr{H}\to \overline{\mathbb R}$ indexed by $\theta\in \Theta$, such that
\begin{equation}\label{linearexpression} \Sbep[X]=\max_{\theta\in \Theta} E_{\theta}[X] \; \text{ for } X \in \mathscr{H} \text{ with } \Sbep[X] \text{ being finite}. \end{equation}
Moreover, for each $X\in \mathscr{H}$, there exists $\theta_X\in \Theta$ such that $\Sbep[X]=E_{\theta_X}[X]$ if $\Sbep[X]$ is finite.

\begin{definition} ({\em See Peng (2008, 2019)})

\begin{description}
  \item[ \rm (i)] ({\em Identical distribution}) Let $\bm X_1$ and $\bm X_2$ be two $n$-dimensional random vectors defined
respectively in sub-linear expectation spaces $(\Omega_1, \mathscr{H}_1, \Sbep_1)$
  and $(\Omega_2, \mathscr{H}_2, \Sbep_2)$. They are called identically distributed, denoted by $\bm X_1\overset{d}= \bm X_2$  if
$$ \Sbep_1[\varphi(\bm X_1)]=\Sbep_2[\varphi(\bm X_2)], \;\; \forall \varphi\in C_{b,Lip}(\mathbb R^n). $$
 A sequence $\{X_n;n\ge 1\}$ of random variables is said to be identically distributed if $X_i\overset{d}= X_1$ for each $i\ge 1$.
\item[\rm (ii)] ({\em Independence})   In a sub-linear expectation space  $(\Omega, \mathscr{H}, \Sbep)$, a random vector $\bm Y =
(Y_1, \ldots, Y_n)$, $Y_i \in \mathscr{H}$ is said to be independent to another random vector $\bm X =
(X_1, \ldots, X_m)$ , $X_i \in \mathscr{H}$ under $\Sbep$  if for each test function $\varphi\in C_{b,Lip}(\mathbb R^m \times \mathbb R^n)$
we have
$ \Sbep [\varphi(\bm X, \bm Y )] = \Sbep \big[\Sbep[\varphi(\bm x, \bm Y )]\big|_{\bm x=\bm X}\big].$
 \item[\rm (iii)] ({\em Independent random variables}) A sequence of random variables $\{X_n; n\ge 1\}$
 is said to be independent, if  $X_{i+1}$ is independent to $(X_{1},\ldots, X_i)$ for each $i\ge 1$.
 \end{description}
\end{definition}

 Let $(\Omega, \mathscr{H}, \Sbep)$ be a sub-linear expectation space.  We denote   $(\Capc,\cCapc)$ to be a pair of  capacities with the properties that
\begin{equation}\label{eq1.3}
   \Sbep[f]\le \Capc(A)\le \Sbep[g]\;\;
\text{ if } f\le I_A\le g, f,g \in \mathscr{H} \text{ and } A\in \mathcal F,
\end{equation}
 $$ \Capc \text{ is sub-additive in sense that } \Capc(A\bigcup B)\le \Capc(A)+\Capc(B)  \text{ for all } A,B\in \mathcal F  $$
and $\cCapc(A):= 1-\Capc(A^c)$, $A\in \mathcal F$.
We call $\Capc$ and $\cCapc$ the upper and the lower capacity, respectively. In general, we can choose $\Capc$ as
\begin{equation}\label{eq1.5} \Capc(A):=\inf\{\Sbep[\xi]: I_A\le \xi, \xi\in\mathscr{H}\},\;\; \forall A\in \mathcal F.
\end{equation}
To distinguish this   capacity from others, we denote it by $\upCapc$, and $\lowCapc(A)=1-\upCapc(A)$. $\upCapc$ is the largest capacity satisfying \eqref{eq1.3}.

When there exists  a family of probability measure on $(\Omega,\mathscr{F})$ such that
\begin{equation}\label{eq1.7} \Sbep[X]=\sup_{P\in \mathscr{P}}P[X]=:\sup_{P\in \mathscr{P}}\int  XdP ,
\end{equation} $\Capc$ can be defined as
\begin{equation}\label{eq1.6} \Capc(A)=\sup_{P\in \mathscr{P}}P(A).
\end{equation}
We denote this capacity by $\Capc^{\mathscr{P}}$, and $\cCapc^{\mathscr{P}}(A)=1-\Capc^{\mathscr{P}}(A)$.

 Also, we define the  Choquet integrals/expecations $(C_{\Capc},C_{\cCapc})$  by
$$ C_V[X]=\int_0^{\infty} V(X\ge t)dt +\int_{-\infty}^0\left[V(X\ge t)-1\right]dt $$
with $V$ being replaced by $\Capc$ and $\cCapc$ respectively.
If $\Capc_1$ on the sub-linear expectation space $(\Omega_1,\mathscr{H}_1,\Sbep_1)$ and $\Capc_2$ on the sub-linear expectation space $(\Omega_2,\mathscr{H}_2,\Sbep_2)$  are two capacities having the property \eqref{eq1.3}, then for any random variables $X_1\in \mathscr{H}_1$ and $\tilde X_2\in \mathscr{H}_2$ with $X_1\overset{d}=\tilde X_2$, we have
\begin{equation}\label{eqV-V}
\Capc_1(X_1\ge x+\epsilon)\le \Capc_2(\tilde X_2\ge x), \;\;  \cCapc_1(X_1\ge x+\epsilon)\le \cCapc_2(\tilde X_2\ge x) \text{ for all } \epsilon>0 \text{ and } x
\end{equation}
and
$C_{\Capc_1}(X_1)=C_{\Capc_2}(\tilde X_2).
$

\begin{definition}\label{def3.1}
   A function $V:\mathcal F\to [0,1]$ is called to be  countably sub-additive if
$$ V\Big(\bigcup_{n=1}^{\infty} A_n\Big)\le \sum_{n=1}^{\infty}V(A_n) \;\; \forall A_n\in \mathcal F. $$
\end{definition}

Since the capacity $\upCapc$ defined as in \eqref{eq1.5}  may be not countably sub-additive, we  consider its countably sub-additive extension.
\begin{definition}   A  countably sub-additive extension $\outCapc$  of $\upCapc$   is defined by
\begin{equation}\label{outcapc} \outCapc(A)=\inf\Big\{\sum_{n=1}^{\infty}\upCapc(A_n): A\subset \bigcup_{n=1}^{\infty}A_n\Big\},\;\; \outcCapc(A)=1-\outCapc(A^c),\;\;\; A\in\mathcal F,
\end{equation}
where $\upCapc$ is defined as in \eqref{eq1.5}.
\end{definition}

As shown in Zhang (2016a), $\outCapc$ is countably sub-additive, $\outCapc(A)\le \upCapc(A)$ and,  if $V$ is also a sub-additive (resp. countably sub-additive) capacity satisfying
\begin{equation}\label{eqVprop} V(A)\le \Sbep[g]  \text{ whenever }  I_A\le g\in \mathscr{H},
\end{equation}
 then $V(A)\le \upCapc$ (resp. $V(A)\le \outCapc(A)$.

 \begin{definition} Another countably sub-additive capacity generated by $\Sbep$ can be defined as follows:
 \begin{equation}\label{tildecapc} \tildeCapc(A)=\inf\Big\{\lim_{n\to\infty}\Sbep[\sum_{i=1}^n g_i]: I_A\le \sum_{n=1}^{\infty}g_n, 0\le g_n\in\mathscr{H}\Big\},\;\;\; A\in\mathcal F.
\end{equation}
\end{definition}

 $\tildeCapc$ is  countably sub-additive and  has the property \eqref{eqVprop}, and so, $\tildeCapc(A)\le \outCapc(A)$. Furthermore,
 if $\Sbep$ has the form \eqref{eq1.7}, then
$$ \Capc^{ \mathscr{P}}(A)=\sup_{P\in \mathscr{P}}P(A)\le \tildeCapc(A)\le \outCapc(A), \;\; A\in\mathcal{F}.  $$

Finally, for real numbers $x$ and $y$, denote $x\vee y=\max(x,y)$, $x\wedge y=\min(x,y)$, $x^+=\max(0,x)$, $x^-=\max(0,-x)$ and $\log x=\ln \max(e,x)$. For a random variable $X$, because $XI\{|X|\le c\}$  may be not in $\mathscr{H}$, we will truncate it in the form $(-c)\vee X\wedge c$ denoted by $X^{(c)}$. We  denote
$$\breve{\mathbb E}[X]=\lim_{c\to \infty} \Sbep[X^{(c)}]$$
if the limit exists.

   \section{The law of the iterated logarithm}\label{sectLIL}
\setcounter{equation}{0}

 Let   $\{Y_n; n\ge 1\}$  be a sequence of independent and identically distributed random variables in the sub-linear expectation space $(\Omega,\mathscr{H},\Sbep)$ with a capacity satisfying \eqref{eq1.3}.
   Set $S_n=\sum_{i=1}^nY_i$.  The following law of the iterated logarithm is obtained by Zhang (2021a).

   \begin{theorem} \label{thLILiid}
 Suppose
  \begin{equation}\label{eqLILmomentcondition1} C_{\Capc}\left[\frac{Y_1^2}{\log\log|Y_1|}\right]<\infty,
    \end{equation}
    \begin{equation}\label{eqLILmeanzerocondition}\breve{\mathbb E}[Y_1]=\breve{\mathbb E}[-Y_1]=0  \end{equation}
    and
    \begin{equation}\label{eqLILmomentcondition2}\overline{\sigma}^2=\lim_{c\to \infty}\Sbep[Y_1^2\wedge c]<\infty. \end{equation}
    Then
\begin{equation}\label{eqthLILiid.3}
\outCapc\left( \limsup_{n\to \infty}\frac{|S_n|}{\sqrt{2n\log\log n}}> \overline{\sigma} \right)=0,
\end{equation}
and
\begin{equation}\label{eqthLILiid.4}\outcCapc\left([-\overline{\sigma}, \;\;\overline{\sigma}]\supset C\Big\{\frac{S_n}{\sqrt{2n\log\log n}}\Big\} \supset
[-\underline{\sigma}, \;\;\underline{\sigma}]\right)=1,
\end{equation}
where $C \{x_n\} $ denotes the cluster set of a sequence of $\{x_n\}$ in $\mathbb R$,  $\underline{\sigma}^2=\lim_{c\to \infty}\cSbep[Y_1^2\wedge c]$.

Further, suppose the following condition:
\begin{description}
      \item[\rm (CC)]      The sub-linear expectation $\Sbep$  on $\mathscr{H}_b$  satisfies
\begin{equation} \label{eqexpressbyP} \Sbep[X]=\sup_{P\in \mathscr{P}}P[X], \; X\in \mathscr{H}_b
\end{equation}
where $\mathscr{H}_b=\{f\in \mathscr{H}; f \text{ is bounded}\}$, $\mathscr{P}$ is a  countable-dimensionally weakly compact family of probability measures on $(\Omega,\sigma(\mathscr{H}))$ in sense that, for any $X_1,X_2,\ldots \in \mathscr{H}_b$ and any  sequence $\{P_n\}\subset \mathscr{P}$ there are a subsequence $\{n_k\}$ and a probability measure $P\in \mathscr{P}$ for which
$$ \lim_{k\to \infty} P_{n_k}[\varphi(X_1,\ldots,X_d)]= P[\varphi(X_1,\ldots,X_d)],\; \varphi\in C_{b,Lip}(\mathbb R^d), d\ge 1.
$$
\end{description}
 Then for $V=\Capc^{\mathscr{P}}$, $\mathbb C^{\ast}$ or $\outCapc$,
 \begin{equation}\label{eqthLILnew2.3}
 V\left(  C \left\{\frac{S_n}{\sqrt{2n\log\log n}}\right\}=[-\overline{\sigma},\overline{\sigma}]\right)=1.
\end{equation}

Conversely, if (CC) is satisfied and for $V=\Capc^{\mathscr{P}}$, $\mathbb C^{\ast}$ or $\outCapc$,
$$V\left(  \limsup_{n\to \infty} \frac{|S_n|}{\sqrt{2n\log\log n}}=\infty\right)<1, $$
then \eqref{eqLILmomentcondition1}-\eqref{eqLILmomentcondition2} hold.
\end{theorem}

Zhang (2021a) also showed that there is a copy $\{\widetilde{Y}_n;n\ge 1\}$ of $\{Y_n;n\ge 1\}$ such that
\begin{equation} \label{eqthLILnew4.4}
\widetilde V\left(   C \left\{\frac{\sum_{i=1}^n \tilde Y_i}{\sqrt{2n\log\log n}}\right\}=[-\sigma, \sigma]\right)=\begin{cases}
  1, & \text{ when } \sigma\in [\underline{\sigma},\overline{\sigma}],\\
  0, & \text{ when } \sigma\not\in [\underline{\sigma},\overline{\sigma}].
  \end{cases}
\end{equation}

The purpose of this note  is to show that \eqref{eqthLILnew4.4} remains true for the original sequence $\{Y_n;n\ge 1\}$.
 \begin{theorem}  Suppose that the condition (CC) is satisfied. Assume \eqref{eqLILmomentcondition1}, \eqref{eqLILmeanzerocondition} and
   \begin{equation}\label{eqLILmomentcondition3} \underline{\sigma}^2>0 \; \text{ and }\;  \overline{\sigma}^2<\infty.
    \end{equation}
     Then for $V=\Capc^{\mathscr{P}}$, $\mathbb C^{\ast}$ or $\outCapc$,
 \begin{equation} \label{eq:LILcompact}
V\left(   C \left\{\frac{\sum_{i=1}^n   Y_i}{\sqrt{2n\log\log n}}\right\}=[-\sigma, \sigma]\right)=\begin{cases}
  1, & \text{ when } \sigma\in [\underline{\sigma},\overline{\sigma}],\\
  0, & \text{ when } \sigma\not\in [\underline{\sigma},\overline{\sigma}].
  \end{cases}
\end{equation}
 \end{theorem}

 \begin{remark} Zhang (2021b) showed that the condition (CC) is equivalent to    that $\{P\bm X^{-1}:P\in \mathscr{P}\}$   is a weakly compact family of probability measures on the metric space $\mathbb R^{\infty}$ for any $\bm X=(X_1,X_2,\ldots)$ with $X_i\in \mathscr{H}_b$, $i=1,2,\ldots$, and it is also equivalent to  that $\Sbep$ is regular on $\mathscr{H}_b$ in sense that $\Sbep[\varphi]\to 0$ whenever  $\mathscr{H}_n\ni \varphi\searrow 0$.
 \end{remark}
  \section{Proofs}\label{sectProof}
  \setcounter{equation}{0}

Recall that   $\{Y_n; n\ge 1\}$  is a sequence of  independent and identically distributed random variables in the sub-linear expectation space $(\Omega,\mathscr{H},\Sbep)$. Let $p>2$. Denote  $t_j=\sqrt{2\log\log j}$ and $b_j=\alpha_j   \sqrt{j}/\sqrt{2\log\log j}$, where $\alpha_j>0$  is  specified such that $\alpha_j\to 0$, $b_j\nearrow \infty$  and $\alpha_j^{1-p}t_j^{-2}\to 0$. Let $V_n^2=\sum_{j=1}^n (Y_j^2\wedge b_j^2)$.

When $\sigma\not\in [\underline{\sigma},\overline{\sigma}]$, \eqref{eq:LILcompact} follows from \eqref{eqthLILiid.3} and \eqref{eqthLILiid.4} immediately.
When $\sigma\in [\underline{\sigma},\overline{\sigma}]$, \eqref{eq:LILcompact} will  follow  from a self-normalized law of the iterated logarithm and a law of large numbers which are stated as Propositions \ref{prop:1} and \ref{prop:2}, respectively.
\begin{proposition} \label{prop:1} Suppose there is a   family of probability measures $\mathscr{P}$ on $(\Omega,\sigma(\mathscr{H}))$ such that the sub-linear expectation $\Sbep$   satisfies  \eqref{eqexpressbyP}. Assume that \eqref{eqLILmomentcondition1}, \eqref{eqLILmeanzerocondition} and \eqref{eqLILmomentcondition3} are satisfied.
Then
\begin{equation}\label{eq:prop1.1} \cCapc^{\mathscr{P}}\left(C\left\{ \frac{\sum_{i=1}^n Y_i}{V_n\sqrt{2\log\log n}}\right\}=[-1,1]\right)=1.
\end{equation}
 \end{proposition}

 \begin{proposition} \label{prop:2}  Suppose  $\overline{\sigma}^2<\infty$.
Then
\begin{equation}\label{eq:prop2.1}\outcCapc\left(\underline{\sigma}^2\le \liminf_{n\to \infty}\frac{V_n^2}{n}\le \limsup_{n\to \infty}\frac{V_n^2}{n}\le \overline{\sigma}^2\right)=1.
\end{equation}
Furthermore, if the sub-linear expectation satisfies the condition (CC), then
\begin{equation}\label{eq:prop2.2} \Capc^{\mathscr{P}}\left(\lim_{n\to \infty}\frac{V_n^2}{n}=\sigma^2\right)=1\;\; \text{for any}\;\; \sigma^2\in [\underline{\sigma}^2,\overline{\sigma}^2].
\end{equation}
 \end{proposition}

Suppose  the condition \eqref{eqLILmomentcondition1} is strengthened to
 \begin{equation}\label{eqLILmomentcondition4} C_{\Capc}\left(Y_1^2\right)<\infty.
    \end{equation}
    By the strong law of large numbers (c.f. Zhang (2021b)), with capacity one under $\outcCapc$ we have
$$ \limsup_{n\to \infty} \frac{\sum_{j=1}^n Y_j^2-V_n^2}{n}\le \limsup_{n\to \infty} \frac{\sum_{i=1}^n (Y_j^2-C)^+}{n}\le \Sbep[(Y_1^2-C)^+]\le  C_{\Capc}[(Y_1^2-C)^+]\to 0
$$
as $C\to \infty$, and so $\sum_{j=1}^n Y_j^2/V_n^2\to 1$. Hence, we have a corollary of Proposition \ref{prop:1}.
 \begin{corollary} Suppose there is a   family of probability measures $\mathscr{P}$ on $(\Omega,\sigma(\mathscr{H}))$ such that the sub-linear expectation $\Sbep$   satisfies  \eqref{eqexpressbyP}. Assume that \eqref{eqLILmomentcondition4}, \eqref{eqLILmeanzerocondition} and \eqref{eqLILmomentcondition3} are satisfied.
 Then
 \begin{equation}\label{eq:corollary1.1} \cCapc^{\mathscr{P}}\left(C\left\{ \frac{\sum_{i=1}^n Y_i}{U_n\sqrt{2\log\log n}}\right\}=[-1,1]\right)=1,
\end{equation}
where  $U_n^2=\sum_{i=1}^n Y_j^2$.
 \end{corollary}

Before proving the propositions, we give some remarks.

\begin{remark}\label{remark3}
Under  general moment conditions, Zhang (2016b) showed the self-normalized law of the iterated logarithm as
\begin{equation}\label{eq:Zhang2016b} \Capc\left(C\left\{ \frac{\sum_{i=1}^n Y_i}{U_n\sqrt{2\log\log n}}\right\}=[-1,1]\right)=1,
\end{equation}
under the assumption that $\Capc$ is continuous. Zhang (2021a) showed that $\Capc$ is not continuous  unless the sub-linear expectation is reduced to a linear expectation. However, with the same arguments as in Zhang (2021a,b), by applying Lemma \ref{lem:BCconverse} (ii) instead of the original Borel-Cantelli Lemma,  \eqref{eq:Zhang2016b} remains true if  $\Capc$ is  replaced by a countably sub-additive capacity $\outCapc$, $\mathbb C^{\ast}$ or $\Capc^{\mathscr{P}}$ when the condition (CC) is satisfied.

 After checking  a special case that the random variables are G-normal, Zhang (2016b) conjectured that  the upper capacity $\Capc$ can be  replaced by the lower capacity $\cCapc$ in  \eqref{eq:Zhang2016b}. The equation \eqref{eq:corollary1.1} gives partly answer to this conjecture. \eqref{eq:corollary1.1} and \eqref{eq:Zhang2016b} are interesting because   the self-normalization  eliminates the variance-uncertainty.
\end{remark}

\begin{remark}\label{remark1} From \eqref{eq:prop1.1} and \eqref{eq:prop2.1}, it follows that
\begin{equation}\label{eq:LILunderlowerC}\cCapc^{\mathscr{P}}\left([-\overline{\sigma}, \;\;\overline{\sigma}]\supset C\Big\{\frac{S_n}{\sqrt{2n\log\log n}}\Big\} \supset
[-\underline{\sigma}, \;\;\underline{\sigma}]\right)=1,
\end{equation}
which is similar to \eqref{eqthLILiid.4}. However, \eqref{eq:LILunderlowerC} is also a directly corollary of \eqref{eqthLILiid.4} since $\outcCapc(A)\le \cCapc^{\mathscr{P}}(A)$. We conjecture that  \eqref{eq:prop1.1} remains true when $\cCapc^{\mathscr{P}}$ is replaced by $\outcCapc$.
\end{remark}
\begin{remark}\label{remark2} Recently, by the means of the martingale
analogue of the Kolmogorov law of the iterated logarithm in the classical probability theory, Gao et al. (2022) established  the laws of the iterated logarithm for general independent  random variables which are not necessarily identically distributed and may have nonzero means. As a result,  \eqref{eq:LILunderlowerC} is established   when    the moment condition \eqref{eqLILmomentcondition1} is replaced by \eqref{eqLILmomentcondition4}.

\end{remark}

 For proving the Propositions \ref{prop:1} and \ref{prop:2}, we need several Lemmas.
The first lemma gives the  exponential  inequalities.
\begin{lemma}\label{lem:ExpIneq}   Suppose that  $\{X_1,\ldots, X_n\}$ is a sequence  of independent random variables
on $(\Omega, \mathscr{H}, \Sbep)$. Set
$A_n(p,y)=\sum_{i=1}^n\Sbep[(X_i^+\wedge y )^p]$ and
$ \breve{B}_{n,y}=\sum_{i=1}^n \vSbep[(X_i\wedge y)^2]$.   Then  for all $p\ge 2$, $x,y>0$, $0<\delta\le 1$,
\begin{align}\label{eq:ExpIneq.1}
& \upCapc\Big( \max_{k\le n}  \sum_{i=1}^k(X_i-\vSbep[X_i])\ge x\Big)
\le    \upCapc\big(\max_{k\le n} X_k> y \big)
+\exp\left\{-\frac{x^2}{2(xy+\breve{B}_{n,y})   }\right\}.
\end{align}
and
\begin{align}\label{eq:ExpIneq.2}
& \upCapc\Big( \max_{k\le n}  \sum_{i=1}^k(X_i-\vSbep[X_i])\ge x\Big) \nonumber\\
\le &  \upCapc\big(\max_{k\le n} X_k> y \big)
  +2\exp\{p^p\}\Big\{\frac{A_n(p,y)}{y^p} \Big\}^{\frac{\delta x}{10y}}
+\exp\left\{-\frac{x^2}{2\breve{B}_{n,y}(1+\delta)   }\right\}.
\end{align}
 \end{lemma}

{\bf Proof}.    The proof  is the same as that of (3.1) and (3.2) of  Zhang (2021a) if we note
$$ \Sbep[e^{t (X_k\wedge y)}]=\vSbep[e^{t (X_k\wedge y)}]\le 1+t\vSbep[X_k]+\frac{e^{ty}-1-t y}{y^2} \vSbep[(X_k\wedge y)^2], \; y>0. \;\;\;\; \Box$$

 The following exponential  inequality for martingales can be found in  Pe\~{n}a (1999). See also Theorem 9.12 of Pe\~{n}a et al. (2009).
\begin{lemma}\label{lem:ExpIneqM} Let $\{X_n; n\ge 1\}$ be a martingale difference sequence in a probability space $(\Omega,\mathcal{F},\pr)$ with respect to the
filtration $\{\mathcal F_n\}$ such that $|X_n|\le  c$ a.s.   Then
$$ \pr\left( \sum_{i=1}^n X_i\ge x \text{ and }  \sum_{i=1}^n \ep[X_i^2|\mathcal{F}_{i-1}]\le y\; \text{ for some } n\right)
\le \exp\left\{ -\frac{x^2}{2(xc+y)}\right\}.   $$
\end{lemma}

The following   self-normalized law of the iterated logarithm for martingales  is  Lemma 13.8 of   Pe\~{n}a et al. (2009).
\begin{lemma}\label{lem:selfLILM} Let $\{X_n; n\ge 1\}$ be a martingale difference sequence in a probability space $(\Omega,\mathcal{F},\pr)$ with respect to the
filtration $\{\mathcal F_n\}$ such that $|X_n|\le  m_n$ a.s. for some $\mathcal F_{n-1}$-measurable random variable
$m_n$, with $U_n\to \infty$  and $m_n/\{U_n(\log\log U_n)^{-1/2}\}\to 0$  a.s., where $U_n^2=\sum_{i=1}^n X_i^2$.  Then
$$ \limsup_{n\to \infty}\frac{\sum_{i=1}^n X_i}{U_n(2\log\log U_n)^{1/2}}=1\;\; a.s. $$
\end{lemma}

The following lemma gives the relation between the sub-linear expectation  and the conditional expectation under a probability, the  proof of which  can be found in  Guo and Li (2021) (see also Hu et al. (2021) and Gao et al. (2022)).
\begin{lemma}\label{lem3} Let $\{X_n;n\ge 1\}$ be a sequence of independent random variables in the sub-linear expectation space $(\Omega,\mathscr{H},\Sbep)$ with \eqref{eqexpressbyP}. Denote
$$ \mathcal{F}_n=\sigma(X_1,\ldots,X_n) \; \text{ and }\; \mathcal{F}_0=\{\emptyset,\Omega). $$
Then for each $P\in\mathscr{P}$, we have
$$ E_P\left[\varphi(X_n)|\mathcal{F}_{n-1}\right]\le \Sbep[\varphi(X_n)]\;\; a.s., \;\; \varphi\in C_{b,Lip}(\mathbb R). $$
\end{lemma}

The next lemma is the Borel-Cantelli lemma for a countably sub-additive capacity,  the proof of which is trivial   and omitted.
\begin{lemma}\label{lem:BCdirect} Let $V$ be a   countably sub-additive capacity and $\sum_{n=1}^{\infty}V(A_n)<\infty$. Then
$$ V(A_n\; i.o.)=0, \;\; \text{ where } \{ A_n\; i.o.\}=\bigcap_{n=1}^{\infty}\bigcup_{i=n}^{\infty}A_i. $$
\end{lemma}

The next lemma is  the converse part of the Borel-Cantelli lemma (see Zhang (2021b)).
\begin{lemma}\label{lem:BCconverse} Let $(\Omega,\mathscr{H},\Sbep)$ be a sub-linear expectation  space satisfying the condition (CC).
Denote $\Capc^{\mathscr{P}}(A)=\max_{P\in \mathscr{P}}P(A)$, $\cCapc^{\mathscr{P}}(A)=1-\Capc^{\mathscr{P}}(A)$, $A\in \sigma(\mathscr{H}). $
Suppose that $\{X_n;n\ge 1\}$ is a sequence of independent random variables in  $(\Omega,\mathscr{H},\Sbep)$.
 \begin{description}
   \item[\rm (i)] If $\sum_{n=1}^{\infty}\cCapc^{\mathscr{P}}(X_n>1)<\infty$,  then
$$\cCapc^{\mathscr{P}}\left( X_n>1 \; i.o.\right)=0.
$$
   \item[\rm (ii)]   If $\sum_{n=1}^{\infty}\Capc^{\mathscr{P}}(X_n\ge 1)=\infty$, then
$$ \Capc^{\mathscr{P}}\left( X_n\ge 1\;\; i.o. \right)=1.
$$
 \end{description}
\end{lemma}

The last lemma is Lemma 6.1 of Zhang (2021a).
\begin{lemma} \label{lem2} Suppose $X\in \mathscr{H}$.
\begin{description}
  \item[\rm (i)]
 For any $\delta>0$,
$$ \sum_{n=1}^{\infty} \Capc\big(|X|\ge \delta \sqrt{n\log\log n} \big)<\infty \;\; \Longleftrightarrow C_{\Capc}\left[\frac{X^2}{\log\log|X|}\right]<\infty.
$$
 \item[\rm (ii)]
   If $C_{\Capc}\left[\frac{X^2}{\log\log|X|}\right]<\infty$, then for any $\delta>0$ and $p>2$,
$$ \sum_{n=1}^{\infty} \frac{\Sbep\big[\big(|X|\wedge (\delta \sqrt{n\log\log n})\big)^p\big]}{(n\log\log n)^{p/2}}<\infty. $$
 \item[\rm (iii)] $C_{\Capc}\left[\frac{X^2}{\log\log|X|}\right]<\infty$, then for any $\delta>0$,
 $$ \Sbep[X^2\wedge (2\delta n\log\log n)]=o(\log \log n) $$
 and $$ \breve{\mathbb E}[(|X|-\delta \sqrt{2  n\log\log n})^+]=o(\sqrt{\log\log n/n}). $$
\end{description}
\end{lemma}

{\bf Proof of Proposition \ref{prop:1}}. Let $P\in \mathscr{P}$, $\mathcal{F}_n=\sigma(Y_1,\ldots,Y_n)$, $\mathcal{F}_0=\{\emptyset,\Omega\}$. By Lemma \ref{lem3},
$$ E_P[Z_j|\mathcal{F}_{j-1}]\le \Sbep[Z_j]=\Sbep[ Y_1^{( b_j)}]\to\vSbep[Y_1]=0. $$
Similarly
$$ E_P[-Z_j|\mathcal{F}_{j-1}]\le \Sbep[-Z_j]=\Sbep[-Y_1^{( b_j)}]\to\vSbep[-Y_1]=0. $$
On the other hand, by \eqref{eq:prop2.1},
\begin{equation}\label{eq:prop:1.1} V_n^2=\sum_{i=1}^n Z_j^2\approx n \;\; a.s.
\end{equation}
It follows that
$$ U_n^2=:\sum_{i=1}^n (Z_i-E_P[Z_j|\mathcal{F}_{i-1}])^2\sim V_n^2 \approx n\;\;a.s. $$
It follows that $b_n/\{U_n(\log\log U_n)^{-1/2}\}\to 0$ a.s.
By Lemma \ref{lem:selfLILM}, we have
\begin{align}\label{eq:LILforZ1} & P\left(\limsup_{n\to\infty}\frac{\sum_{i=1}^n (Z_i-E_P[Z_j|\mathcal{F}_{i-1}])}{V_n(2\log\log n)^{1/2}}=1\right)\\
=& P\left(\limsup_{n\to\infty}\frac{\sum_{i=1}^n (Z_i-E_P[Z_j|\mathcal{F}_{i-1}])}{U_n(2\log\log U_n)^{1/2}}=1\right)=1. \nonumber
\end{align}

Next, we will show that $Z_i-E_P[Z_j|\mathcal{F}_{i-1}])$ in \eqref{eq:LILforZ1} can be replaced by $Y_j$. Denote $ d_n=\sqrt{2n\log\log n}$.
Let  $\lambda > 1$. Denote $n_k=[\lambda^k]$ and $I(k)=\{n_k+1,\ldots,n_{k+1}\}$. Then
 $n_k/n_{k+1}\to  {1}/{\lambda} $, $d_{n_k}/d_{n_{k+1}}\to 1/\sqrt{\lambda}$.

Note $p>2$. By Lemma \ref{lem2} (ii), we have
\begin{equation} \label{eqproofthLIL5.7}  \sum_{k=1}^{\infty} \frac{   \Lambda_{n_k,n_{k+1}}(p)}{d_{n_{k+1}}^p}<\infty,
\end{equation}
where
$$\Lambda_{n_k,n_{k+1}}(p)=\sum_{j\in I(k)}\Sbep[\big((|Y_j|\wedge d_{n_{k+1}}\big)^p]. $$
 Let
\begin{equation} \label{eqproofthLIL5.8}   \mathbb N_1=\left\{k\in \mathbb N; \frac{  \Lambda_{n_k, n_{k+1}}(p)}{d_{n_{k+1}}^p}\le t_{n_{k+1}}^{-2p}\right\}.
\end{equation}
We have for $k\in \mathbb N_1$,
\begin{align}
 &\frac{\sum_{j\in I(k)}
\Sbep[ |Y_j^{(d_{n_{k+1}})}-Z_j|]}{d_{n_{k+1}}}
\le   C\frac{\Lambda_{n_k,n_{k+1}}(p)}{d_{n_{k+1}}^p}  \alpha_{n_{k+1}}^{1-p}t_{n_{k+1}}^{2p-2} \le  C\alpha_{n_{k+1}}^{1-p}t_{n_{k+1}}^{-2}\to 0,\label{eqproofthLIL5.11}\\
 &\frac{\sum_{j\in I(k)}
  \Sbep[(Y_j^{(d_{n_{k+1}})}-Z_j)^2] }{ n_{k+1} }
   \le   C\frac{\Lambda_{n_k,n_{k+1}}(p)}{d_{n_{k+1}}^p}  \alpha_{n_{k+1}}^{2-p}t_{n_{k+1}}^{2p-2} \le   C \alpha_{n_{k+1}}^{2-p}t_{n_{k+1}}^{-2}\to 0,\label{eqproofthLIL5.12}
\end{align}
by noting $\alpha_j\to 0$ such that $\alpha_j^{1-p}t_j^{-2}\to 0$. Define
$$ Z_{j,1}=Z_j, \text{ if } j\in I(k)\; \text{ with }\;  k\in \mathbb N_1, \text{ and } 0  \; \text{ otherwise}, $$
$$ Z_{j,2}=0, \text{ if } j\in I(k)\; \text{ with }\;  k\in \mathbb N_1, \text{ and } Z_j \; \text{ otherwise}. $$
Then $Z_j=Z_{j,1}+Z_{j,2}$. Note $|Z_j-E_P[Z_j|\mathcal F_{j-1}|\le 2 b_j\le 2b_{n_{k+1}}$, $j\in I(k)$, and
$$ E_P[(Z_j-E_P[Z_j|\mathcal F_{j-1})^2|\mathcal F_{j-1}]\le  E_P[Z_j^2|\mathcal F_{j-1}]
\le  \Sbep[Z_j^2]\le  \overline{\sigma}^2 \; a.s. $$
by Lemma \ref{lem3}. By Lemma \ref{lem:ExpIneqM} with $c=2b_{n_{k+1}}$, $y=(n_{k+1}-n_k)\overline{\sigma}^2$ and $x=\epsilon d_{n_{k+1}}$, we have
\begin{align*}
& P\left(\max_{n\in I(k)}|\sum_{i=n_k+1}^n (Z_i-E_P[Z_i|\mathcal{F}_{i-1}])| \ge \epsilon d_{n_{k+1}}\right)\\
\le & 2\exp\left\{ -\frac{\epsilon^2 d_{n_{k+1}}^2}{2(\epsilon d_{n_{k+1}}\cdot 2 b_{n_{k+1}}+(n_{k+1}-n_k)\overline{\sigma}^2)}\right\}
\le 2\exp\left\{-\epsilon^{\prime} t_{n_{k+1}}^2\right\}.
\end{align*}
It follows that
\begin{align*}
& \sum_{k=1}^{\infty} P\left(\max_{n\in I(k)}|\sum_{i=n_k+1}^n (Z_{i,2}-E_P[Z_{i,2}|\mathcal{F}_{i-1}])| \ge \epsilon d_{n_{k+1}}\right)\\
\le &\sum_{k\not\in \mathbb N_1}  2\exp\left\{-\epsilon^{\prime} t_{n_{k+1}}^2\right\}\le \sum_{k\not\in \mathbb N_1}  C t_{n_{k+1}}^{-2p}<\infty,
\end{align*}
by \eqref{eqproofthLIL5.7} and \eqref{eqproofthLIL5.8}.  Hence, by the Borel-Cantelli lemma,
$$ \frac{\max_{n\in I(k)}|\sum_{i=n_k+1}^n (Z_{i,2}-E_P[Z_{i,2}|\mathcal{F}_{i-1}])| }{d_{n_{k+1}}}\to 0\;\; a.s. \text{ under } P, $$
which implies
$$ \frac{ \sum_{i=1}^n (Z_{i,2}-E_P[Z_{i,2}|\mathcal{F}_{i-1}]) }{d_n}\to 0\;\; a.s. \text{ under } P. $$
Hence, by \eqref{eq:prop:1.1},
\begin{equation}\label{eq:LILforZ2} \frac{ \sum_{i=1}^n (Z_{i,2}-E_P[Z_{i,2}|\mathcal{F}_{i-1}]) }{V_n(2\log\log n)^{1/2}}\to 0\;\; a.s.  \text{ under } P.
\end{equation}
By Lemma \ref{lem3} and noticing that $\vSbep[Y_j]=0$, we have
\begin{align*} E_P[Z_j|\mathcal{F}_{j-1}]\le  & \Sbep[Z_j]\le \vSbep[Y_j]+\vSbep[Y_j^{(d_{n_{k+1}})}-Y_j]+\Sbep[Z_j-Y_j^{(d_{n_{k+1}})}] \\
\le & \Sbep[|Y_j^{(d_{n_{k+1}})}-Z_j|]+\vSbep[(|Y_1|-d_{n_{k+1}})^+].
\end{align*}
Similarly,
$$E_P[-Z_j|\mathcal{F}_{j-1}]\le   \Sbep[|Y_j^{(d_{n_{k+1}})}-Z_j|]+\vSbep[(|Y_1|-d_{n_{k+1}})^+].$$
By \eqref{eqproofthLIL5.11} and Lemma \ref{lem2} (iii),
$$ \sum_{j\in I(k)}|E_P[Z_j|\mathcal{F}_{j-1}]|= o(d_{n_{k+1}}), \;\; k\in \mathbb N_1\;\; a.s. \text{ under } P,
 $$
 which implies that
 $$ \frac{\sum_{i=1}^n |E_P[Z_{j,1} |\mathcal{F}_{j-1}]|}{d_n}\to 0\;\; a.s. \text{ under } P. $$
 Hence, by \eqref{eq:prop:1.1},
\begin{equation}\label{eq:LILforZ3} \frac{ \sum_{i=1}^n |E_P[Z_{j,1} |\mathcal{F}_{j-1}]| }{V_n(2\log\log n)^{1/2}}\to 0\;\; a.s. \text{ under } P.
\end{equation}
Combing \eqref{eq:LILforZ1}, \eqref{eq:LILforZ2} and \eqref{eq:LILforZ3} yields that
\begin{equation}\label{eq:LILforZ}   P\left(\limsup_{n\to\infty}\frac{\sum_{i=1}^n Z_{i,1}}{V_n(2\log\log n)^{1/2}}=1\right)=1.
\end{equation}

Now, we consider $Y_j$. Let $x=\epsilon a_{n_{k+1}}$ and $y=\epsilon^{\prime}a_{n_{k+1}}$, where  $\epsilon^{\prime}>0$ is chosen  such that  $x/(10 y)\ge   1 $. By    the inequality \eqref{eq:ExpIneq.2} in Lemma \ref{lem:ExpIneq},
\begin{align*}
 & \Capc\left(\max_{n\in I(k)}\sum_{j=n_k+1}^n Y_j  \ge \epsilon d_{n_{k+1}}\right)  \\
\le & \exp\left\{-\frac{\epsilon^2}{2(1+1)} t_{n_{k+1}}^2\right\}+\sum_{j\in I(k)} \Capc\big(|Y_j|\ge \epsilon^{\prime} d_{n_{k+1}}\big)
 +C \frac{\sum_{j\in I(k)}\Sbep[|Y_j|^p\wedge d_{n_{k+1}}^{p}]}{d_{n_{k+1}}^{p}}\\
 \le &  \sum_{j\in I(k)} \Capc\big(|Y_1|\ge \epsilon^{\prime} d_j\big)
 +C \frac{\Lambda_{n_k,n_{k+1}}(p)}{d_{n_{k+1}}^{p}}, \;\; k\not\in \mathbb N_1.
\end{align*}
 It follows that
\begin{equation}\label{eq:LILforY1} \sum_{k\not\in \mathbb N_1}\Capc\left(\max_{n\in I(k)}|\sum_{j=n_k+1}^n Y_j|  \ge \epsilon d_{n_{k+1}}\right)<\infty,
\end{equation}
by Lemma \ref{lem2} and \eqref{eqproofthLIL5.7}.
For $k\in \mathbb N_1$, note \eqref{eqproofthLIL5.11} and \eqref{eqproofthLIL5.12}. By the inequality \eqref{eq:ExpIneq.2}  of  Lemma \ref{lem:ExpIneq} with $\delta=1$, it follows   that for $k$ large enough,
\begin{align*}
&\Capc\Big(\sum_{j\in I(k)}|Y_j-Z_j|\ge \epsilon d_{n_{k+1}}\Big) \\
\le & \Capc\Big(\sum_{j\in I(k)}|  Y_j^{( d_{n_{k+1}})}-Z_j|\ge \epsilon d_{n_{k+1}}\Big)+\sum_{j\in I(k)}\Capc(|Y_j|>d_{n_{k+1}})  \\
\le & \Capc\Big(\sum_{j\in I(k)}\big(|Y_j^{( d_{n_{k+1}})}-Z_j|-\Sbep[|Y_j^{( d_{n_{k+1}})}-Z_j|]\big)\ge \epsilon d_{n_{k+1}}/2\Big)
  +\sum_{j\in I(k)}\Capc(|Y_j|>d_{n_{k+1}})  \\
  \le &  \exp\left\{-\frac{(\epsilon d_{n_{k+1}}/2)^2}{4 \sum_{j\in I(k)}\vSbep[(Y_j^{( d_{n_{k+1}})}-Z_j)^2]} \right\}
    +C \frac{\Lambda_{n_k,n_{k+1}}(p)}{d_{n_{k+1}}^p} +2\sum_{j\in I(k)}\Capc(|Y_j|>\epsilon^{\prime} d_j)\\
\le &  \exp\left\{-2 t_{n_{k+1}}^2\right\}
    +C \frac{\Lambda_{n_k,n_{k+1}}(p)}{d_{n_{k+1}}^p} +2\sum_{j\in I(k)}\Capc(|Y_j|>\epsilon^{\prime} d_j).
\end{align*}
It follows that
\begin{equation}\label{eq:LILforY2}
\sum_{k\in \mathbb N_1 } \Capc\Big(\sum_{j\in I(k)}|Y_j-Z_j|\ge \epsilon d_{n_{k+1}}\Big)<\infty.
\end{equation}
Combining \eqref{eq:LILforY1} and \eqref{eq:LILforY2} yields that
\begin{equation}\label{eq:LILforY} \sum_{k=1}^{\infty}\Capc\left(\max_{n\in I(k)}|\sum_{j=n_k+1}^n (Y_j-Z_{j,1})|  \ge \epsilon d_{n_{k+1}}\right)<\infty.
\end{equation}
Note $\max_{n\in I(k)}|\sum_{j=n_k+1}^n (Y_j-Z_{j,1})|\in \mathscr{H}$. Hence, by \eqref{eqV-V},
$$ \sum_{k=1}^{\infty}\upCapc\left(\max_{n\in I(k)}|\sum_{j=n_k+1}^n (Y_j-Z_{j,1})|  \ge \epsilon d_{n_{k+1}}\right)<\infty,
$$
which implies
$$ \outcCapc\left(\frac{\max_{n\in I(k)}|\sum_{j=n_k+1}^n (Y_j-Z_{j,1})| }{ d_{n_{k+1}}}\to 0\right)=1, $$
by the countably sub-additivity of $\outCapc$, $\outCapc(A)\le \upCapc(A)$ and the Borel-Cantelli lemma-Lemma \ref{lem:BCdirect}.  Hence
$$ \outcCapc\left( \frac{\sum_{i=1}^n (Y_i-Z_{i,1})  }{ d_n}\to 0\right)=1,  $$
which, together with \eqref{eq:prop2.1}, yields that
\begin{equation}\label{eq:LILforY-Z}
\outcCapc\left( \frac{\sum_{i=1}^n (Y_i-Z_{i,1})  }{ V_n(2\log\log n)^{1/2}}\to 0\right)=1.
\end{equation}
From \eqref{eq:LILforZ} and \eqref{eq:LILforY-Z}, it fellows that
$$ P\left(\limsup_{n\to\infty}\frac{\sum_{i=1}^n Y_i}{V_n(2\log\log n)^{1/2}}=1\right)=1.
$$
Replacing $Y_i$ by $-Y_i$ in the above equality yields
$$ P\left(\liminf_{n\to\infty}\frac{\sum_{i=1}^n Y_i}{V_n(2\log\log n)^{1/2}}=-1\right)=1.
$$
Also, it is easily shown that
$$\outcCapc\left( \lim_{n\to \infty} \frac{\max_{i\le n}|Y_i| }{ V_n(2\log\log n)^{1/2}}=\lim_{n\to \infty} \frac{\max_{i\le n}|Y_i| }{ d_n}=0\right)=1.
$$
It follows that
$$ P\left(C\left\{\frac{\sum_{i=1}^n Y_i}{V_n(2\log\log n)^{1/2}}\right\}=[-1,1]\right)=1\;\;\text{ for all } P\in \mathscr{P}.
$$
The proof of Proposition \ref{prop:1} is completed by noting $\cCapc^{\mathscr{P}}(A)=\inf_{P\in \mathscr{P}}P(A)$.
$\Box$

{\bf Proof of Proposition \ref{prop:2}}.  Let $n_k=2^k$.  Applying the inequality \eqref{eq:ExpIneq.1} in Lemma \ref{lem:ExpIneq} with $x=\epsilon n_k$, $y=2b_{n_k}^2$ and $X_j=Z_j^2$ yields
\begin{align}\label{eq:prop2:proof1}
& \upCapc\left(\max_{n\le n_k}\sum_{i=1}^n (Z_i^2-\Sbep[Z_i^2])\ge \epsilon n_k\right)
\le \exp\left\{-\frac{\epsilon^2 n_k^2}{2(2\epsilon n_kb_{n_k}+\sum_{i=1}^{n_k} \Sbep[Z_j^4])}\right\}\nonumber\\
\le & \exp\left\{-\frac{\epsilon^2 n_k^2}{2(2\epsilon n_kb_{n_k}^2+n_k b_{n_k}^2 \overline{\sigma}^2})\right\}\le \exp\{-2\log\log n_k\}.
\end{align}
Hence
$$ \sum_k \upCapc\left(\max_{n\le n_k}\sum_{i=1}^n (Z_i^2-\Sbep[Z_i^2])\ge \epsilon n_k\right)<\infty \text{ for all } \epsilon>0. $$
It follows that
$$ \outcCapc\left(\limsup_{k\to \infty} \frac{\max_{n\le n_k}\sum_{i=1}^n (Z_i^2-\Sbep[Z_i^2])}{n_k}\le 0\right)=1 $$
by the countable sub-additivity of $\outCapc$, $\outCapc(A)\le \upCapc(A)$ and the Borel-Cantelli lemma-Lemma \ref{lem:BCdirect}. The above equality  implies
$$ \outcCapc\left(\limsup_{n\to \infty} \frac{ \sum_{i=1}^n (Z_i^2-\Sbep[Z_i^2])}{n}\le 0\right)=1. $$
Note $\Sbep[Z_j^2]=\Sbep[Y_1^2\wedge b_j]\to \vSbep[Y_1^2]=\overline{\sigma}^2$. We have
$$ \outcCapc\left(\limsup_{n\to \infty} \frac{ \sum_{i=1}^n Z_i^2}{n}\le \overline{\sigma}^2\right)=1. $$
Replacing $Z_j^2$ by $-Z_j^2$ in the above equality yields
$$ \outcCapc\left(\liminf_{n\to \infty} \frac{ \sum_{i=1}^n Z_i^2}{n}\ge \underline{\sigma}^2\right)=1. $$
 The proof of \eqref{eq:prop2.1} is completed.

 For \eqref{eq:prop2.2},  we use the arguments in Zhang (2021b). By \eqref{linearexpression},  there are finite additive linear expectations $E_{j,1}$ and $E_{j,2}$ such that
 $$ E_{j,1}[Z_j^2]=\cSbep[Z_j^2] \text{ and }  E_{j,2}[Z_j^2]=\Sbep[Z_j^2]. $$
 Note $\cSbep[Z_j^2]\to \underline{\sigma}^2$ and $\Sbep[Z_j]\to \overline{\sigma}^2$. For $\sigma^2\in [\underline{\sigma}^2,\overline{\sigma}^2]$, there exists $0\le\alpha_j\le 1$ such that $\alpha_j\cSbep[Z_j^2]+(1-\alpha_j)\Sbep[Z_j^2]\to \sigma^2$. Let $E_j=\alpha_jE_{j,1}+(1-\alpha_j)E_{j,2}$. Then $E_j$ is a finite additive linear expectation with $E_j\le \Sbep$. By Proposition 2.1 of Zhang (2021b),   we can find a new sub-linear space $(\widetilde{\Omega},\widetilde{\mathscr{H}},\widetilde{\mathbb E})$ defined on a metric space $\widetilde{\Omega}=\mathbb R^{\infty}$, with a  copy $\{\widetilde Y_n;n\ge 1\}$ on $(\widetilde{\Omega},\widetilde{\mathscr{H}},\widetilde{\mathbb E})$ of $\{Y_n;n\ge 1\}$  and a probability measure $Q$ on $\widetilde{\Omega}$ such that
   $\{\widetilde Y_n;n\ge 1\}$ is a sequence of independent random variables under $Q$,
\begin{equation}\label{eq:copy1}
E_Q\left[\varphi(\widetilde Y_i)\right]=E_i\left[\varphi(Y_i)\right]\; \text{ for all } \varphi\in C_{b,Lip}(\mathbb R),
\end{equation}
$$
E_Q\left[\varphi(\widetilde Y_1,\ldots,\widetilde Y_d)\right]\le \widetilde{\mathbb E}\left[\varphi(\widetilde Y_1,\ldots,\widetilde Y_d)\right]=\Sbep\left[\varphi(Y_1,\ldots,Y_d)\right]\; \text{ for all } \varphi\in C_{b,Lip}(\mathbb R^d)
$$
and
\begin{equation}\label{eq:copy3}
\widetilde{v}(B)\le Q(B)\le \widetilde{V}(B) \; \text{ for all } B\in\sigma(\widetilde Y_1,\widetilde Y_2,\ldots),
\end{equation}
where
$$ \widetilde{V}(A)=\widetilde{\Capc}^{\widetilde{\mathscr{P}}}(A)=\sup_{P\in\widetilde{\mathscr{P}}}P(A)\; \text{ and } \widetilde{v}(A)=1-\widetilde{V}(A), \;\;  A\in \widetilde{\mathcal F},$$
and $\widetilde{\mathscr{P}}$ is the family of all probability measures $P$ on $(\widetilde{\Omega},\widetilde{\mathcal F})$ with the property
$$ E_P[\varphi]\le \widetilde{\mathbb E}[\varphi] \; \text{ for bounded }\varphi\in  \widetilde{\mathscr{H}}.  $$
By \eqref{eq:copy1}, $E_Q[\widetilde Y_j^4\wedge b_j^4]=E_j[  Y_j^4\wedge b_j^4]\le b_j^2 E_j[Y_j^2\wedge b_j^2]\sim  b_j^2\sigma^2$. Similarly to \eqref{eq:prop2:proof1},
 \begin{align*}
& Q\left(\max_{n_k+1\le n\le  n_{k+1}}\Big|\sum_{j=n_k+1}^n (\widetilde Y_j^2\wedge b_j^2-E_Q[\widetilde Y_j^2\wedge b_j^2])\Big|\ge \epsilon n_{k+1}\right) \\
\le &2\exp\left\{-\frac{\epsilon^2 n_{k+1}^2}{2(2\epsilon n_{k+1}b_{n_{k+1}}^2+\sum_{i=n_k+1}^{n_{k+1}}E_Q[\widetilde Y_j^4\wedge b_j^4])}\right\}\nonumber\\
\le &2 \exp\left\{-\frac{\epsilon^2 n_{k+1}^2}{2(2\epsilon n_{k+1}b_{n_{k+1}}^2+2n_k b_{n_{k+1}}^2  \sigma^2)}\right\}\le 2\exp\{-2\log\log n_k\}
\end{align*}
for $k$ large enough. Note $E_Q[\widetilde Y_j^2\wedge b_j^2]=E_j[ Y_j^2\wedge b_j^2] \to \sigma^2$. It follows that
$$\sum_{k=1}^{\infty} Q\left(\max_{n_k+1\le n\le  n_{k+1}}\Big|\sum_{j=n_k+1}^n (\widetilde Y_j^2\wedge b_j^2-\sigma^2)\Big|\ge 2\epsilon n_{k+1}\right) <\infty
\; \text{ for all } \epsilon>0. $$
Then, there exists a sequence $0<\epsilon_k\searrow 0$ such that
 $$\sum_{k=1}^{\infty} Q\left(\max_{n_k+1\le n\le  n_{k+1}}\Big|\sum_{j=n_k+1}^n (\widetilde Y_j^2\wedge b_j^2-\sigma^2)\Big|\ge  \epsilon_k n_{k+1}\right) <\infty.  $$
By \eqref{eq:copy3}, we have
$$\sum_{k=1}^{\infty} \widetilde{v}\left(\max_{n_k+1\le n\le  n_{k+1}}\Big|\sum_{j=n_k+1}^n (\widetilde Y_j^2\wedge b_j^2-\sigma^2)\Big|\ge  \epsilon_k n_{k+1}\right) <\infty.  $$
Note $(Y_1,\ldots,Y_{n_{k+1}})\overset{d}=(\widetilde Y_1,\ldots, \widetilde Y_{n_{k+1}})$. By \eqref{eqV-V}, we have
  $$\sum_{k=1}^{\infty} \cCapc^{\mathscr{P}}\left(\max_{n_k+1\le n\le  n_{k+1}}\Big|\sum_{j=n_k+1}^n (  Y_j^2\wedge b_j^2-\sigma^2)\Big|>  2\epsilon_k n_{k+1}\right) <\infty.  $$
Note the independence. By  the Borel-Cantelli lemma-Lemma \ref{lem:BCconverse},
$$\cCapc^{\mathscr{P}}\left(A_k\; i.o.  \right)=0 \; \text{ with } \; A_k=\left\{\frac{\max\limits_{n_k+1\le n\le  n_{k+1}}\Big|\sum_{j=n_k+1}^n (  Y_j^2\wedge b_j^2-\sigma^2)\Big|}{n_{k+1}}>  2\epsilon_k \right\}. $$
  Notice that on the event $(A_k\; i.o.)^c$,
  $$ \lim_{k\to \infty}\frac{\max\limits_{n_k+1\le n\le  n_{k+1}}\Big|\sum_{j=n_k+1}^n (  Y_j^2\wedge b_j^2-\sigma^2)\Big|}{n_k}=0, $$
  which implies $ \lim_{n\to \infty} \frac{\sum_{j=1}^n (  Y_j^2\wedge b_j^2)}{n}=\sigma^2$. \eqref{eq:prop2.2} is proved. $\Box$

\bigskip


\end{document}